**Fourier series representations of the logarithms of the Euler gamma function and the Barnes multiple gamma functions**

Donal F. Connon

dconnon@btopenworld.com

25 March 2009

**Abstract**

Kummer's Fourier series for $\log \Gamma(t)$ is well known, having been discovered in 1847. In this paper we develop a corresponding Fourier series for logarithm of the Barnes double gamma function (and the method may be easily extended to the higher order multiple gamma functions). Some applications of these Fourier series are explored.

## 1. Introduction

We recall the Hasse identity for the Hurwitz zeta function [27] which holds for all $s \in \mathbf{C}$ except $s = 1$

$$(1.1) \qquad \varsigma(s,t) = \frac{1}{s-1} \sum_{n=0}^{\infty} \frac{1}{n+1} \sum_{k=0}^{n} \binom{n}{k} \frac{(-1)^k}{(t+k)^{s-1}}$$

and with $s \rightarrow 1-s$ this becomes

$$(1.2) \qquad s\varsigma(1-s,t) = -\sum_{n=0}^{\infty} \frac{1}{n+1} \sum_{k=0}^{n} \binom{n}{k} (-1)^k (t+k)^s$$

We have the well-known Hurwitz's formula for the Fourier expansion of the Hurwitz zeta function $\varsigma(s,t)$ as reported in Titchmarsh's treatise [38, p.37]

$$(1.3) \qquad \varsigma(s,t) = 2\Gamma(1-s)\left[ \sin\left(\frac{\pi s}{2}\right) \sum_{n=1}^{\infty} \frac{\cos 2n\pi t}{(2\pi n)^{1-s}} + \cos\left(\frac{\pi s}{2}\right) \sum_{n=1}^{\infty} \frac{\sin 2n\pi t}{(2\pi n)^{1-s}} \right]$$

where $\mathrm{Re}(s) < 0$ and $0 < t \leq 1$. In 2000, Boudjelkha [12] showed that this formula also applies in the region $\mathrm{Re}(s) < 1$. It may be noted that when $t = 1$ this reduces to Riemann's functional equation for $\varsigma(s)$. Letting $s \rightarrow 1-s$ we may write (1.3) as

$$(1.4) \qquad \varsigma(1-s,t) = 2\Gamma(s)\left[ \cos\left(\frac{\pi s}{2}\right) \sum_{n=1}^{\infty} \frac{\cos 2n\pi t}{(2\pi n)^{s}} + \sin\left(\frac{\pi s}{2}\right) \sum_{n=1}^{\infty} \frac{\sin 2n\pi t}{(2\pi n)^{s}} \right]$$

$$= 2\Gamma(s)\sum_{n=1}^{\infty}\frac{\cos[\pi s/2 - 2n\pi t]}{(2\pi n)^s}$$

which is valid for

$$(\sigma < 0,\ 0 < t \le 1;\ 0 < \sigma,\ t < 1)$$

The derivation of (1.3) has been simplified by Zhang and Williams [41].

## 2. Kummer's Fourier series representation of the gamma function

Multiplying (1.4) by $s$, we see that

$$(2.1) \qquad f(s,t) = s\varsigma(1-s,t) = 2\Gamma(s+1)\sum_{n=1}^{\infty}\frac{\cos[\pi s/2 - 2n\pi t]}{(2\pi n)^s}$$

$$(2.2) \qquad = -\sum_{n=0}^{\infty}\frac{1}{n+1}\sum_{k=0}^{n}\binom{n}{k}(-1)^k (t+k)^s$$

Differentiation of (2.2) results in

$$(2.3) \qquad \frac{\partial^p}{\partial s^p}f(s,t) = -\sum_{n=0}^{\infty}\frac{1}{n+1}\sum_{k=0}^{n}\binom{n}{k}(-1)^k (t+k)^s \log^p(t+k)$$

and we have the particular value at $s = 1$

$$(2.4) \qquad \frac{\partial^p}{\partial s^p}f(s,t)\bigg|_{s=1} = -\sum_{n=0}^{\infty}\frac{1}{n+1}\sum_{k=0}^{n}\binom{n}{k}(-1)^k (t+k) \log^p(t+k)$$

We also have from (2.1)

$$(2.5)$$
$$\frac{\partial}{\partial s}f(s,t) = -2\Gamma(s+1)\sum_{n=1}^{\infty}\frac{(\pi/2)\sin[\pi s/2 - 2n\pi t]}{(2\pi n)^s} - 2\Gamma(s+1)\sum_{n=1}^{\infty}\frac{\cos[\pi s/2 - 2n\pi t]\log(2\pi n)}{(2\pi n)^s}$$

$$+ 2\Gamma'(s+1)\sum_{n=1}^{\infty}\frac{\cos[\pi s/2 - 2n\pi t]}{(2\pi n)^s}$$

and thus

$$\frac{\partial}{\partial s}f(s,t)\bigg|_{s=1} = -\frac{1}{2}\sum_{n=1}^{\infty}\frac{\cos 2n\pi t}{n} - \sum_{n=1}^{\infty}\frac{\sin 2n\pi t \log(2\pi n)}{\pi n} + \Gamma'(2)\sum_{n=1}^{\infty}\frac{\sin 2n\pi t}{\pi n}$$

We note the familiar trigonometric series shown in Carslaw's book [13, p.241] (elementary derivations are also contained in [18])

$$(2.6) \qquad -\log(2\sin\pi t) = \sum_{n=1}^{\infty} \frac{\cos 2n\pi t}{n} \qquad (0 < t < 1)$$

$$(2.7) \qquad \pi\left(\frac{1}{2} - t\right) = \sum_{n=1}^{\infty} \frac{\sin 2n\pi t}{n} \qquad (0 < t < 1)$$

Using these identities results in

$$\frac{\partial}{\partial s} f(s,t)\bigg|_{s=1} = \frac{1}{2}\log(2\sin\pi t) + [1 - \gamma - \log(2\pi)]\left(\frac{1}{2} - t\right) - \sum_{n=1}^{\infty} \frac{\sin 2n\pi t \log(n)}{\pi n}$$

and we then have

$$-\sum_{n=0}^{\infty} \frac{1}{n+1} \sum_{k=0}^{n} \binom{n}{k} (-1)^k (t+k)\log(t+k)$$

$$= \frac{1}{2}\log(2\sin\pi t) + [1 - \gamma - \log(2\pi)]\left(\frac{1}{2} - t\right) - \sum_{n=1}^{\infty} \frac{\sin 2n\pi t \log(n)}{\pi n}$$

We showed in [14] that

$$(2.8) \qquad \log\Gamma(t) = \sum_{n=0}^{\infty} \frac{1}{n+1} \sum_{k=0}^{n} \binom{n}{k} (-1)^k (t+k)\log(t+k) + \frac{1}{2} - t + \frac{1}{2}\log(2\pi)$$

and we therefore obtain Kummer's Fourier series [9] for the log gamma function (which, because we relied on (2.6) and (2.7), is only valid for $0 < t < 1$)

$$(2.9) \qquad \log\Gamma(t) = \frac{1}{2}\log\frac{\pi}{\sin\pi t} + [\gamma + \log(2\pi)]\left(\frac{1}{2} - t\right) + \frac{1}{\pi}\sum_{n=1}^{\infty} \frac{\log n}{n}\sin 2\pi n t$$

Reference to (2.6) and (2.7) confirms that (2.9) is properly described as a Fourier series expansion for $\log\Gamma(t)$. In 1985 Berndt [9] gave an elementary proof of this Fourier series expansion, which was originally derived by Kummer in 1847.

The series immediately gives rise to the familiar result

$$(2.10) \qquad \log\Gamma\left(\frac{1}{2}\right) = \frac{1}{2}\log\pi$$

It may be noted that it is not possible to differentiate (2.9) because, as is easily seen, the resulting infinite series is divergent. Notwithstanding this, a trigonometric expansion (which is not a Fourier series) exists for the digamma function $\psi(t)$ as shown in (8.7) below.

Alternatively, differentiating (1.1) gives us

$$(2.11) \qquad (s-1)\varsigma'(s,t) + \varsigma(s,t) = -\sum_{n=0}^{\infty} \frac{1}{n+1} \sum_{k=0}^{n} \binom{n}{k} (-1)^k \frac{\log(t+k)}{(t+k)^{s-1}}$$

and evaluation at $s = 0$ produces

$$(2.12) \qquad \varsigma'(0,t) = \varsigma(0,t) + \sum_{n=0}^{\infty} \frac{1}{n+1} \sum_{k=0}^{n} \binom{n}{k} (-1)^k (t+k) \log(t+k)$$

We have the well known relationship between the Hurwitz zeta function and the Bernoulli polynomials $B_n(t)$ (for example, see Apostol's book [5, pp. 264-266])

$$(2.13) \qquad \varsigma(-m,t) = -\frac{B_{m+1}(t)}{m+1} \quad \text{for } m \in \mathbf{N}_o$$

and it may be noted that this identity may also be deduced from (1.2) because we have shown in [19]

$$(2.14) \qquad B_{m+1}(t) = \sum_{n=0}^{\infty} \frac{1}{n+1} \sum_{k=0}^{n} \binom{n}{k} (-1)^k (t+k)^{m+1}$$

In particular, from (2.13) we have

$$\varsigma(0,t) = -B_1(t) = \frac{1}{2} - t$$

From (2.12) we then see that

$$(2.15) \qquad \varsigma'(0,t) = \frac{1}{2} - t + \sum_{n=0}^{\infty} \frac{1}{n+1} \sum_{k=0}^{n} \binom{n}{k} (-1)^k (t+k) \log(t+k)$$

and, comparing this with (2.8), we have therefore deduced Lerch's identity [9]

$$(2.16) \qquad \varsigma'(0,t) = \log \Gamma(t) - \frac{1}{2} \log(2\pi)$$

Since $\varsigma'(0,1) = \varsigma'(0) = -\frac{1}{2}\log(2\pi)$ this may be expressed as

$$\varsigma'(0,t) - \varsigma'(0) = \log\Gamma(t)$$

Lerch established the above relationship between the gamma function and the Hurwitz zeta function in 1894 (other derivations are contained in, for example, Berndt's paper [9] and [14]).

As noted by Berndt [9] we have with $t \rightarrow 1-t$ in (2.9)

$$(2.17) \qquad \log\Gamma(1-t) = \frac{1}{2}\log\frac{\pi}{\sin\pi t} - [\gamma + \log(2\pi)]\left(\frac{1}{2} - t\right) + \frac{1}{\pi}\sum_{n=1}^{\infty}\frac{\log n}{n}\sin 2\pi nt$$

and thus adding (2.9) and (2.17) together we see that

$$\log\Gamma(t) + \log\Gamma(1-t) = \log\frac{\pi}{\sin\pi t}$$

which is simply Euler's reflection formula for the gamma function

$$(2.18) \qquad \Gamma(t)\Gamma(1-t) = \frac{\pi}{\sin\pi t}$$

$\square$

Differentiating (2.8) results in

$$\psi(t) = \sum_{n=0}^{\infty}\frac{1}{n+1}\sum_{k=0}^{n}\binom{n}{k}(-1)^k\log(t+k) + \sum_{n=0}^{\infty}\frac{1}{n+1}\sum_{k=0}^{n}\binom{n}{k}(-1)^k - 1$$

and, since $\sum_{k=0}^{n}\binom{n}{k}(-1)^k = \delta_{n,0}$, we see that [14]

$$(2.19) \qquad \psi(t) = \sum_{n=0}^{\infty}\frac{1}{n+1}\sum_{k=0}^{n}\binom{n}{k}(-1)^k\log(t+k)$$

This result was also recently obtained in a different way by Guillera and Sondow [24].

Differentiating (2.19) gives us

$$(2.20) \qquad \psi'(t) = \sum_{n=0}^{\infty}\frac{1}{n+1}\sum_{k=0}^{n}\binom{n}{k}\frac{(-1)^k}{(t+k)}$$

and, as seen from (1.1), this is equal to $\varsigma(2,t)$.

$$\square$$

Since $\sin(n\pi/2) = -\sin(3n\pi/2)$, letting $t = 1/4$ and $t = 3/4$ respectively in (2.9) and adding the two equations together, we obtain

$$\log\Gamma\left(\frac{1}{4}\right) + \log\Gamma\left(\frac{3}{4}\right) = \log\pi + \frac{1}{2}\log 2$$

which of course may also be easily obtained from Euler's reflection formula (2.18) for the gamma function (or, alternatively, from Legendre's duplication formula for the gamma function [36, p.7]).

Noting that $\sin(n\pi/2) = \sin(5n\pi/2)$ unfortunately does not assist us because $t = 5/4$ falls outside of the region of validity of Kummer's formula (2.9).

## 3. An application of Parseval's theorem

Applying Parseval's theorem [8, p.338] to the Fourier series (2.9) we have

$$(3.1) \quad \int_0^1 \left[ \log\Gamma(t) - \frac{1}{2}\log\frac{\pi}{\sin\pi t} - [\gamma + \log(2\pi)]\left(\frac{1}{2} - t\right) \right]^2 dt = \frac{1}{2\pi^2}\sum_{n=1}^{\infty}\frac{\log^2 n}{n^2} = \frac{\varsigma''(2)}{2\pi^2}$$

or equivalently we have the six component integrals

$$\int_0^1 \log^2\Gamma(t)\,dt + \frac{1}{4}\int_0^1 \log^2\frac{\pi}{\sin\pi t}\,dt + [\gamma + \log(2\pi)]^2 \int_0^1 \left(\frac{1}{2} - t\right)^2 dt$$

$$-\int_0^1 \log\Gamma(t)\log\frac{\pi}{\sin\pi t}\,dt + [\gamma + \log(2\pi)]\int_0^1 \left(\frac{1}{2} - t\right)\log\frac{\pi}{\sin\pi t}\,dt$$

$$-2[\gamma + \log(2\pi)]\int_0^1 \left(\frac{1}{2} - t\right)\log\Gamma(t)\,dt = \frac{\varsigma''(2)}{2\pi^2}$$

In order to determine $\int_0^1 \log^2\Gamma(t)\,dt$ we now evaluate the last five of these integrals in turn.

**Second integral**

We have

$$\int_0^1 \log^2 \frac{\pi}{\sin \pi t}\, dt = \int_0^1 (\log \pi - \log \sin \pi t)^2\, dt$$

$$= \log^2 \pi - 2 \log \pi \int_0^1 \log \sin \pi t\, dt + \int_0^1 \log^2 \sin \pi t\, dt$$

The Log-Sine integrals $\mathrm{Ls}_n(\theta)$ are defined for $n \geq 2$ by

$$(3.2) \qquad \mathrm{Ls}_n(\theta) = -\int_0^\theta \left( \log \left| 2 \sin \frac{t}{2} \right| \right)^{n-1} dt$$

and these integrals have been considered by many authors, including Beumer [10], Lewin [31], Boros and Moll [11, p.245] and Srivastava and Choi [36, p.118]. From [31] we have for example

$$(3.3) \qquad \mathrm{Ls}_2(\pi) = -\int_0^\pi \log \left[ 2 \sin \frac{t}{2} \right] dt = 0$$

$$(3.4) \qquad \mathrm{Ls}_3(\pi) = -\int_0^\pi \log^2 \left[ 2 \sin \frac{t}{2} \right] dt = -\frac{\pi^3}{12}$$

With the substitution $\pi x = \frac{t}{2}$ we see that

$$\int_0^\pi \log \left[ 2 \sin \frac{t}{2} \right] dt = 2\pi \int_0^{1/2} \log \left[ 2 \sin \pi x \right] dx$$

and with the substitution $x = 1 - y$ we have

$$\int_0^{1/2} \log \left[ 2 \sin \pi x \right] dx = \int_{1/2}^1 \log \left[ 2 \sin \pi y \right] dy$$

and hence

$$\int_0^1 \log[2\sin\pi x]\,dx = 2\int_0^{1/2} \log[2\sin\pi x]\,dx = 0$$

This gives us the well-known integral (see also (5.10) below)

(3.5)     $$\int_0^1 \log\sin\pi t\,dt = -\log 2$$

Similarly we have

$$\int_0^\pi \log^2\left[2\sin\frac{t}{2}\right]dt = 2\pi\int_0^{1/2}\log^2[2\sin\pi x]\,dx = \pi\int_0^1 \log^2[2\sin\pi x]\,dx$$

and we then see from (3.4) that

(3.6)     $$\int_0^1 \log^2[2\sin\pi t]\,dt = \frac{\pi^2}{12}$$

We have

$$\int_0^1 \log^2[2\sin\pi t]\,dt = \int_0^1 \log^2\sin\pi t\,dt + 2\log 2\int_0^1 \log\sin\pi t\,dt + \log^2 2$$

and therefore using (3.6) we have

(3.7)     $$\int_0^1 \log^2\sin\pi t\,dt = \frac{\pi^2}{12} + \log^2 2$$

This is the problem published by Bremekamp [10] in 1957.

Alternatively, we could also apply Parseval's theorem to (2.6) and obtain

$$\int_0^1 \log^2(2\sin\pi t)\,dt = \frac{1}{2}\sum_{n=1}^\infty \frac{1}{n^2} = \frac{\pi^2}{12}$$

Applying Parseval's theorem to (2.7) results in

$$\pi^2\int_0^1\left(\frac{1}{2}-t\right)^2 dt = \frac{1}{2}\sum_{n=1}^\infty \frac{1}{n^2}$$

It is easily seen that

$$\int_0^1 \left( \frac{1}{2} - t \right)^2 dt = \frac{1}{12}$$

and we therefore obtain Euler's formula for $\varsigma(2)$

$$\sum_{n=1}^{\infty} \frac{1}{n^2} = \frac{\pi^2}{6}$$

This method could also be applied to the fifth and sixth integrals.

To conclude this part we have

$$\int_0^1 \log^2 \frac{\pi}{\sin \pi t} dt = \int_0^1 (\log \pi - \log \sin \pi t)^2 dt$$

$$= \log^2 \pi - 2\log \pi \int_0^1 \log \sin \pi t \, dt + \int_0^1 \log^2 \sin \pi t \, dt$$

$$= \log^2 \pi + 2\log \pi \log 2 + \frac{\pi^2}{12} + \log^2 2$$

$$= \log^2(2\pi) + \frac{\pi^2}{12}$$

**Third integral**

The third integral is rather basic but, as a generalisation, we note [5, p.276]

$$\int_0^1 B_n^2(t) dt = (-1)^{n+1} \frac{(n!)^2}{(2n)!} B_{2n}$$

and thus

$$\int_0^1 \left( \frac{1}{2} - t \right)^2 dt = \frac{1}{12}$$

**Fourth integral**

We showed in equation (6.123) of [18] that

(3.8) $\qquad \int\limits_0^1 \log \Gamma(x+1) \log \left[ 2\sin(\pi x) \right] dx = \frac{1}{2\pi} \sum\limits_{n=1}^\infty \frac{si(2n\pi)}{n^2} = \frac{1}{2\pi} \sum\limits_{n=1}^\infty \frac{Si(2n\pi)}{n^2} - \frac{1}{4}\varsigma(2)$

where $Si(x)$ is the sine integral function defined by [23, p.878] and [1, p.231] as

$$Si(x) = \int\limits_0^x \frac{\sin t}{t} dt \qquad , \quad Si(0) = 0$$

We have the well-known integral from Fourier series analysis

$$\frac{\pi}{2} = \int\limits_0^\infty \frac{\sin t}{t} dt$$

and therefore defining

$$si(x) = Si(x) - \frac{\pi}{2}$$

we have

$$si(x) = \int\limits_0^x \frac{\sin t}{t} dt - \int\limits_0^\infty \frac{\sin t}{t} dt = -\int\limits_x^\infty \frac{\sin t}{t} dt$$

In equation (6.117j) of [18] it was also shown that

(3.9) $\qquad \frac{1}{2\pi^2} \sum\limits_{n=1}^\infty \frac{Si(2n\pi)}{n^2} = \log A - \frac{1}{4}$

where $A$ is the Glaisher-Kinkelin constant

$$\log A = \frac{1}{12} - \varsigma'(-1)$$

Therefore we have

(3.10) $\qquad \int\limits_0^1 \log \Gamma(x+1) \log \left[ 2\sin(\pi x) \right] dx = \pi \log A - \frac{\pi}{4} - \frac{1}{4}\varsigma(2)$

It is easily seen that

$$\int\limits_0^1 \log \Gamma(x+1) \log \left[ 2\sin(\pi x) \right] dx = \int\limits_0^1 \left( \log x + \log \Gamma(x) \right)\left( \log 2 + \log \sin(\pi x) \right) dx$$

$$= \log 2 \int_0^1 \log x \, dx + \log 2 \int_0^1 \log \Gamma(x) \, dx + \int_0^1 \log x \log \sin(\pi x) \, dx + \int_0^1 \log \Gamma(x) \log \sin(\pi x) \, dx$$

$$= -\log 2 + \frac{1}{2} \log 2 \log(2\pi) + \int_0^1 \log x \log \sin(\pi x) \, dx + \int_0^1 \log \Gamma(x) \log \sin(\pi x) \, dx$$

where we have used Raabe's integral

$$\int_0^1 \log \Gamma(x) \, dx = \frac{1}{2} \log(2\pi)$$

which may also be obtained directly from Alexeiewsky's theorem (4.8) below.

Using (2.6) we have

$$\int_0^1 \log x \log \sin(\pi x) \, dx = -\log 2 \int_0^1 \log x \, dx - \int_0^1 \log x \sum_{n=1}^{\infty} \frac{\cos 2n\pi x}{n} \, dx$$

$$= \log 2 - \sum_{n=1}^{\infty} \frac{1}{n} \int_0^1 \log x \cos 2n\pi x \, dx$$

where we have assumed that it is valid to interchange the order of integration and summation.

Let us now consider the integral

$$\int_0^u \log x . \cos ax \, dx = \log x \frac{\sin ax}{a} \Big|_0^u - \int_0^u \frac{\sin ax}{ax} \, dx$$

We have

$$\lim_{x \to 0} \log x \sin ax = \lim_{x \to 0} \left[ ax \log x \frac{\sin ax}{ax} \right] = 0$$

Therefore we get

$$\int_0^u \log x . \cos ax \, dx = \frac{\sin au \log u}{a} - \int_0^u \frac{\sin ax}{ax} \, dx$$

$$= \frac{\sin au \log u}{a} - \frac{1}{a}\int\limits_{0}^{au}\frac{\sin x}{x}\,dx$$

and hence we get

$$\int\limits_{0}^{u}\log x.\cos ax\,dx = \frac{\sin au \log u}{a} - \frac{Si(au)}{a}$$

Particular cases are as follows

$$\int\limits_{0}^{u}\log x.\cos 2n\pi x\,dx = \frac{\sin 2n\pi u \log u}{a} - \frac{Si(2n\pi u)}{2n\pi}$$

$$\int\limits_{0}^{1}\log x.\cos 2n\pi x\,dx = -\frac{Si(2n\pi)}{2n\pi}$$

Hence we have

$$(3.11) \qquad \int\limits_{0}^{1}\log x \log \sin(\pi x)\,dx = \log 2 + \frac{1}{2\pi}\sum_{n=1}^{\infty}\frac{Si(2n\pi)}{n^{2}}$$

$$= \pi \log A - \frac{\pi}{4} + \log 2$$

Therefore we obtain

$$(3.12) \qquad \int\limits_{0}^{1}\log \Gamma(x) \log \sin(\pi x)\,dx = -\frac{1}{2}\log 2 \log(2\pi) - \frac{\pi^{2}}{24}$$

which was previously determined by Espinosa and Moll [21] in a very different manner. We could also have employed the generalised Parseval's theorem [8, p.343]

$$\frac{1}{\pi}\int\limits_{-\pi}^{\pi} f(x)g(x)dx = \frac{1}{2}a_{0}\alpha_{0} + \sum_{n=1}^{\infty}(a_{n}\alpha_{n} + b_{n}\beta_{n})$$

to evaluate (3.12) by utilising the known Fourier series for the two components of the integrand.

**Fifth integral**

We see that

$$\int_0^1 \left( \frac{1}{2} - t \right) \log \frac{\pi}{\sin \pi t} \, dt = \log \pi \int_0^1 \left( \frac{1}{2} - t \right) dt + \int_0^1 B_1(t) \log \sin \pi t \, dt$$

and this is a particular case of an integral noted by Espinosa and Moll [21]

$$(3.13) \qquad \int_0^1 B_{2n+1}(t) \log \sin \pi t \, dt = 0$$

$$(3.14) \qquad \int_0^1 B_{2n}(t) \log \sin \pi t \, dt = \frac{(-1)^n (2n)! \, \varsigma(2n+1)}{(2\pi)^{2n}}$$

Very elementary proofs of the above integrals are given in [18] where we used the basic identity

$$(3.15) \qquad \int_a^b p(x) \cot(\alpha x / 2) \, dx = 2 \sum_{n=1}^{\infty} \int_a^b p(x) \sin \alpha n x \, dx$$

which, as shown in [18], is valid for a wide class of suitably behaved functions. Specifically we require that $p(x)$ is a twice continuously differentiable function. It should be noted that in the above formula we require either (i) both $\sin(x/2)$ and $\cos(x/2)$ have no zero in $[a, b]$ or (ii) if either $\sin(a/2)$ or $\cos(a/2)$ is equal to zero then $p(a)$ must also be zero. Condition (i) is equivalent to the requirement that $\sin x$ has no zero in $[a, b]$.

We then have

$$\int_0^1 \left( t - \frac{1}{2} \right) \log \sin \pi t \, dt = 0$$

In fact this can be shown much more directly by using the substitution $x = t - 1/2$ and then noting that the integrand of the resulting integral is an odd function.

This gives us the fifth integral

$$\int_0^1 \left( \frac{1}{2} - t \right) \log \frac{\pi}{\sin \pi t} \, dt = 0$$

**Sixth integral**

We note that

$$\int_0^1 \left( \frac{1}{2} - t \right) \log \Gamma(t)\, dt = -\int_0^1 B_1(t) \log \Gamma(t)\, dt$$

which is a particular case of an integral noted by Espinosa and Moll [21]

(3.16)    $$\int_0^1 B_{2n-1}(t) \log \Gamma(t)\, dt = \frac{B_{2n}}{2n} \left[ \frac{\varsigma'(2n)}{\varsigma(2n)} - \log(2\pi) - \gamma \right]$$

(3.17)    $$\int_0^1 B_{2n}(t) \log \Gamma(t)\, dt = \frac{(-1)^{n+1}(2n)!\,\varsigma(2n+1)}{2(2\pi)^{2n}} = -\varsigma'(-2n)$$

and we therefore obtain

$$\int_0^1 \left( \frac{1}{2} - t \right) \log \Gamma(t)\, dt = -\frac{1}{12} \left[ \frac{\varsigma'(2)}{\varsigma(2)} - \log(2\pi) - \gamma \right]$$

Using (4.5) we have

$$\frac{\varsigma'(2)}{\varsigma(2)} = \log(2\pi) + \gamma - 1 + 12\varsigma'(-1)$$

and thus

(3.18)    $$\int_0^1 \left( \frac{1}{2} - t \right) \log \Gamma(t)\, dt = \frac{1}{12} - \varsigma'(-1) = \log A$$

The above collection of the above five integrals now enables us to evaluate the first integral.

**First integral**

Using the above we determine that

(3.19)    $$\int_0^1 \log^2 \Gamma(t)\, dt = \frac{\gamma^2}{12} + \frac{\pi^2}{48} + \frac{1}{6} \gamma \log(2\pi) + \frac{1}{3} \log^2(2\pi) - [\gamma + \log(2\pi)] \frac{\varsigma'(2)}{\pi^2} + \frac{\varsigma''(2)}{2\pi^2}$$

which Espinosa and Moll [21] also showed, admittedly with much less effort.

## 4. Fourier series representation of the Barnes double gamma functions

We have from (2.5) in the case where $s = 2$

$$\frac{\partial}{\partial s} f(s,t)\Big|_{s=2} = -\frac{1}{2\pi} \sum_{n=1}^{\infty} \frac{\sin 2n\pi t}{n^2} + \frac{1}{\pi^2} \log(2\pi) \sum_{n=1}^{\infty} \frac{\cos 2n\pi t}{n^2}$$

$$+ \frac{1}{\pi^2} \sum_{n=1}^{\infty} \frac{\cos 2n\pi t \log n}{n^2} - \frac{1}{\pi^2} \left(\frac{3}{2} - \gamma\right) \sum_{n=1}^{\infty} \frac{\cos 2n\pi t}{n^2}$$

From (2.3) we note that

$$\frac{\partial}{\partial s} f(s,t)\Big|_{s=2} = -\sum_{n=0}^{\infty} \frac{1}{n+1} \sum_{k=0}^{n} \binom{n}{k} (-1)^k (t+k)^2 \log(t+k)$$

The Barnes double gamma function $\Gamma_2(x) = 1/G(x)$ defined, inter alia, by [36, p.25]

$$(4.1) \quad G(1+t) = (2\pi)^{1/2} \exp\left[-\frac{1}{2}(\gamma t^2 + t^2 + t)\right] \prod_{k=1}^{\infty} \left\{ \left(1 + \frac{t}{k}\right)^k \exp\left(\frac{t^2}{2k} - t\right) \right\}$$

and it is easily seen that $G(1) = 1$.

It was also shown in [14] that the Barnes double gamma function could be expressed as the logarithmic series

$$(4.2) \quad \log G(1+t) = -\frac{1}{2} \sum_{n=0}^{\infty} \frac{1}{n+1} \sum_{k=0}^{n} \binom{n}{k} (-1)^k (t+k)^2 \log(t+k) + t \log \Gamma(t) + \frac{1}{4} B_2(t) + \varsigma'(-1)$$

where $B_n(t)$ are the Bernoulli polynomials.

We therefore have the trigonometric series

$$\log G(1+t) = -\frac{1}{4\pi} \sum_{n=1}^{\infty} \frac{\sin 2n\pi t}{n^2} + \frac{1}{2\pi^2} \left( \log(2\pi) + \gamma - \frac{3}{2} \right) \sum_{n=1}^{\infty} \frac{\cos 2n\pi t}{n^2}$$

$$(4.3)$$

$$+ \frac{1}{2\pi^2} \sum_{n=1}^{\infty} \frac{\cos 2n\pi t \log n}{n^2} + t \log \Gamma(t) + \frac{1}{4} B_2(t) + \varsigma'(-1)$$

which may be expressed as a Fourier series by noting that [4, p.338] (see also (7.9a) and (7.9b) below)

$$B_2(t) = \frac{1}{\pi^2} \sum_{n=1}^{\infty} \frac{\cos 2n\pi t}{n^2}$$

since

$$(4.4a) \qquad B_{2N}(t) = (-1)^{N+1} 2(2N)! \sum_{n=1}^{\infty} \frac{\cos 2n\pi t}{(2\pi n)^{2N}} \qquad , N = 1, 2, \ldots$$

$$(4.4b) \qquad B_{2N+1}(t) = (-1)^{N+1} 2(2N+1)! \sum_{n=1}^{\infty} \frac{\sin 2n\pi t}{(2\pi n)^{2N+1}} \qquad , N = 0, 1, 2, \ldots$$

To obtain a pure Fourier series for $\log G(1+t)$ it would also be necessary to determine the Fourier series expansion for $t \log \Gamma(t)$ using Kummer's identity (2.9).

With $t = 1$ in (4.3) we have since $G(2) = G(1)\Gamma(1) = 1$

$$(4.5) \qquad \frac{1}{2\pi^2} \varsigma'(2) = \frac{1}{12}\big(\log(2\pi) + \gamma - 1\big) + \varsigma'(-1)$$

which may also be easily derived by differentiating the functional equation for the Riemann zeta function (see for example [29]).

Since $\lim_{t \to 0}[t \log \Gamma(t)] = \lim_{t \to 0}[t \log \Gamma(1+t) - t \log t] = 0$, it may be noted that equation (4.3) also applies when $t = 0$ and this also results in (4.5).

<div align="right">□</div>

Letting $t = 1/2$ in (4.3) gives us

$$\log G(3/2) = \frac{1}{2\pi^2}\left(\log(2\pi) + \gamma - \frac{3}{2}\right)\sum_{n=1}^{\infty} \frac{(-1)^n}{n^2} + \frac{1}{2\pi^2}\sum_{n=1}^{\infty} \frac{(-1)^n \log n}{n^2} + \frac{1}{4}\log \pi - \frac{1}{48} + \varsigma'(-1)$$

The alternating Riemann zeta function is defined by

$$\varsigma_a(s) = \sum_{n=1}^{\infty} \frac{(-1)^{n+1}}{n^s}$$

and it is easily seen that

$$\varsigma(s) = \sum_{n=1}^{\infty} \frac{1}{n^s} = \frac{1}{1-2^{-s}}\sum_{n=1}^{\infty} \frac{1}{(2n-1)^s} = \frac{1}{1-2^{-s}}\sum_{n=0}^{\infty} \frac{1}{(2n+1)^s} \qquad , (\mathrm{Re}(s) > 1)$$

$$= \frac{1}{1-2^{1-s}}\sum_{n=1}^{\infty} \frac{(-1)^{n+1}}{n^s} = \frac{1}{1-2^{1-s}}\varsigma_a(s) \qquad , (\mathrm{Re}(s) > 0; \; s \neq 1)$$

We then have $\sum_{n=1}^{\infty} \frac{(-1)^n}{n^2} = -\frac{1}{2}\varsigma(2) = -\frac{\pi^2}{12}$

Differentiating

$$\varsigma_a(s) = (1 - 2^{1-s})\varsigma(s)$$

gives us

$$\varsigma_a'(s) = (1 - 2^{1-s})\varsigma'(s) + 2^{1-s}\varsigma(s)\log 2$$

and thus

$$\varsigma_a'(2) = \sum_{n=1}^{\infty} \frac{(-1)^n \log n}{n^2} = \frac{1}{2}\varsigma'(2) + \frac{1}{2}\varsigma(2)\log 2$$

We then obtain

$$\log G(3/2) = -\frac{1}{24}\left(\log(2\pi) + \gamma - \frac{3}{2}\right) + \frac{1}{4\pi^2}\varsigma'(2) + \frac{1}{24}\log 2 + \frac{1}{4}\log \pi - \frac{1}{48} + \varsigma'(-1)$$

Using (4.5) this becomes

(4.6) $$\log G(3/2) = \frac{1}{24}\log 2 + \frac{1}{4}\log \pi + \frac{3}{2}\varsigma'(-1)$$

Since [36, p.25]

$$G(1+t) = G(t)\Gamma(t)$$

we see that [36, p.26]

(4.7) $$\log G(1/2) = \frac{1}{24}\log 2 - \frac{1}{4}\log \pi + \frac{3}{2}\varsigma'(-1)$$

as originally determined by Barnes [7] in 1899.

□

Using (2.6) and integrating (2.9) results in

$$\int_0^x \log\Gamma(t)\,dt = \frac{x}{2}\log\pi + \frac{1}{4\pi}\sum_{n=1}^{\infty}\frac{\sin 2n\pi x}{n^2} + \frac{1}{2}x\log 2 - \frac{1}{2}[\gamma + \log(2\pi)][B_2(x) - B_2(0)]$$

$$- \frac{1}{2\pi^2}\sum_{n=1}^{\infty}\frac{\log n}{n^2}\cos 2\pi nx + \frac{1}{2\pi^2}\sum_{n=1}^{\infty}\frac{\log n}{n^2}$$

and comparing this with (4.3) we get

$$\int_0^x \log \Gamma(t)dt = \frac{x}{2}\log(2\pi) - \frac{1}{2}[\gamma + \log(2\pi)][B_2(x) - B_2(0)]$$

$$- \log G(1+x) + \frac{1}{2\pi^2}\left(\log(2\pi) + \gamma - \frac{3}{2}\right)\sum_{n=1}^{\infty}\frac{\cos 2n\pi x}{n^2}$$

$$+ x\log\Gamma(x) + \frac{1}{4}B_2(x) + \varsigma'(-1) + \frac{1}{2\pi^2}\sum_{n=1}^{\infty}\frac{\log n}{n^2}$$

Using (4.5) and some algebra, we obtain Alexeiewsky's theorem [36, p.32], a further derivation of which is contained in equation (4.3.85) of [14]

$$(4.8) \qquad \int_0^x \log\Gamma(t)dt = \frac{1}{2}x(1-x) + \frac{x}{2}\log(2\pi) - \log G(1+x) + x\log\Gamma(x)$$

In a similar way, one could integrate (4.2) to obtain the integral $\int_0^x \log G(1+t)dt$.

## 5. The Gosper/Vardi functional equation

With $s = -1$ in (2.11) we obtain

$$(5.1) \qquad \varsigma'(-1,t) = \frac{1}{2}\varsigma(-1,t) + \frac{1}{2}\sum_{n=0}^{\infty}\frac{1}{n+1}\sum_{k=0}^{n}\binom{n}{k}(-1)^k(t+k)^2\log(t+k)$$

Therefore, using (2.17) we see that

$$\varsigma'(-1,t) = -\frac{1}{4}B_2(t) + \frac{1}{2}\sum_{n=0}^{\infty}\frac{1}{n+1}\sum_{k=0}^{n}\binom{n}{k}(-1)^k(t+k)^2\log(t+k)$$

and substituting (4.2) we obtain

$$(5.2) \qquad \log G(1+t) - t\log\Gamma(t) = \varsigma'(-1) - \varsigma'(-1,t)$$

This functional equation was derived by Vardi in 1988 and also by Gosper in 1997 (see Adamchik's paper [3]).

Letting $t \rightarrow 1-t$ in (5.2) gives us

$$\log G(2-t) - (1-t)\log\Gamma(1-t) = \varsigma'(-1) - \varsigma'(-1,1-t)$$

Noting that $\log G(2-t) = \log G(1-t) + \log \Gamma(1-t)$ we obtain

$$(5.3) \qquad \log G(1-t) + t \log \Gamma(1-t) = \varsigma'(-1) - \varsigma'(-1, 1-t)$$

Letting $t \to 1-t$ in (4.3) gives us

$$\log G(2-t) = \frac{1}{4\pi} \sum_{n=1}^{\infty} \frac{\sin 2n\pi t}{n^2} + \frac{1}{2\pi^2} \left( \log(2\pi) + \gamma - \frac{3}{2} \right) \sum_{n=1}^{\infty} \frac{\cos 2n\pi t}{n^2}$$

$$+ \frac{1}{2\pi^2} \sum_{n=1}^{\infty} \frac{\cos 2n\pi t \log n}{n^2} + (1-t) \log \Gamma(1-t) + \frac{1}{4} B_2(1-t) + \varsigma'(-1)$$

and hence we have

$$\log \frac{G(1+t)}{G(1-t)} = -\frac{1}{2\pi} \sum_{n=1}^{\infty} \frac{\sin 2n\pi t}{n^2} + t \log[\Gamma(t)\Gamma(1-t)] + \frac{1}{4}[B_2(t) - B_2(1-t)]$$

Using the well-known property of the Bernoulli polynomials [36, p.60]

$$B_n(1-t) = (-1)^n B_n(t)$$

we then have

$$(5.4) \qquad \log \frac{G(1+t)}{G(1-t)} = -\frac{1}{2\pi} \sum_{n=1}^{\infty} \frac{\sin 2n\pi t}{n^2} + t \log[\Gamma(t)\Gamma(1-t)]$$

This may be written as

$$\log G(1+t) - t \log \Gamma(t) - [\log G(1-t) + t \log \Gamma(1-t)] = -\frac{1}{2\pi} \sum_{n=1}^{\infty} \frac{\sin 2n\pi t}{n^2}$$

and using the Gosper/Vardi identities (5.2) and (5.4) we see that

$$(5.5) \qquad \varsigma'(-1,t) - \varsigma'(-1,1-t) = \frac{1}{2\pi} \sum_{n=1}^{\infty} \frac{\sin 2n\pi t}{n^2}$$

as previously noted by Adamchik [2].

Using (2.6) we may write (5.4) as

$$\log \frac{G(1+t)}{G(1-t)} = \int_0^t \log(2\sin \pi x)\,dx + t \log \pi - t \log \sin(\pi t)$$

which gives us

$$(5.6) \qquad \log \frac{G(1+t)}{G(1-t)} = -t \log \left[ \frac{\sin \pi t}{2\pi} \right] + \int_0^t \log(\sin \pi x)\, dx$$

and using integration by parts, we see that this is equivalent to the following integral formula originally found by Kinkelin [36, p.30] in 1860

$$(5.7) \qquad \log \frac{G(1+t)}{G(1-t)} = t \log(2\pi) - \int_0^t \pi x \cot \pi x\, dx$$

(which is recorded as an exercise in Whittaker and Watson [39, p.264]).

This may also be written for $0 \le t < 1$ as

$$(5.8) \qquad \int_0^t \pi x \cot \pi x\, dx = \varsigma'(-1,t) - \varsigma'(-1,1-t) + t \log(2 \sin \pi t)$$

which was also derived in equation (4.3.158) in [14] by a different method. Integration by parts results in

$$\int_0^t \pi x \cot \pi x\, dx = t \log \sin \pi t - \int_0^t \log \sin \pi x\, dx$$

and we obtain

$$(5.9) \qquad \int_0^t \log(2 \sin \pi x)\, dx = -[\varsigma'(-1,t) - \varsigma'(-1,1-t)]$$

With $t = 1/2$ in (5.9) we rediscover Euler's integral

$$(5.10) \qquad \int_0^{\pi/2} \log \sin x\, dx = -\frac{\pi}{2} \log 2$$

$\square$

From (4.3) and (4.4a) we have

$$\log G(1+t) - t \log \Gamma(t) - \varsigma'(-1) = -\frac{1}{4\pi} \sum_{n=1}^{\infty} \frac{\sin 2n\pi t}{n^2} + \frac{1}{2\pi^2} \big(\log(2\pi) + \gamma - 1\big) \sum_{n=1}^{\infty} \frac{\cos 2n\pi t}{n^2}$$

$$+ \frac{1}{2\pi^2} \sum_{n=1}^{\infty} \frac{\cos 2n\pi t \log n}{n^2}$$

and (5.2) therefore gives us the Fourier series for $\varsigma'(-1, t)$

(5.11)

$$\varsigma'(-1, t) = \frac{1}{4\pi} \sum_{n=1}^{\infty} \frac{\sin 2n\pi t}{n^2} - \frac{1}{2\pi^2} \big(\log(2\pi) + \gamma - 1\big) \sum_{n=1}^{\infty} \frac{\cos 2n\pi t}{n^2} - \frac{1}{2\pi^2} \sum_{n=1}^{\infty} \frac{\cos 2n\pi t \log n}{n^2}$$

which was derived in a different way in equation (4.4.229i) of [17]. This may of course also be obtained more directly just by differentiating (1.3).

Reference should also be made to the 2002 paper by Koyama and Kurokawa, "Kummer's formula for the multiple gamma functions" [29] where they show by a different method that (which are valid for $0 < t < 1$)

(5.12)    $$\log \Gamma_2^*(x) = -\frac{1}{2\pi^2} \sum_{n=1}^{\infty} \frac{\log n}{n^2} \cos 2\pi nx - \frac{\log(2\pi) + \gamma - 1}{2\pi^2} \sum_{n=1}^{\infty} \frac{\cos 2\pi nx}{n^2}$$

$$+ \frac{1}{4\pi} \sum_{n=1}^{\infty} \frac{\sin 2\pi nx}{n^2} + (1-x) \log \Gamma_1(x)$$

(5.13)    $$\log \Gamma_3^*(x) = -\frac{1}{4\pi^3} \sum_{n=1}^{\infty} \frac{\log n}{n^3} \sin 2\pi nx - \frac{2\log(2\pi) + 2\gamma - 3}{8\pi^3} \sum_{n=1}^{\infty} \frac{\sin 2\pi nx}{n^3}$$

$$+ \frac{1}{8\pi^2} \sum_{n=1}^{\infty} \frac{\cos 2\pi nx}{n^3} + \left(\frac{3}{2} - x\right) \log \Gamma_2^*(x) - \frac{1}{2}(1-x)^2 \log \Gamma_1^*(x)$$

where $\Gamma_1^*(x) = \dfrac{\Gamma(x)}{\sqrt{2\pi}}$. It should however be noted that the multiple gamma functions

$\Gamma_n^*(x)$ considered by Koyama and Kurokawa [29] are not the same as those traditionally employed by Barnes [7], Adamchik [3] etc.

## 6. Fourier series for the alternating Hurwitz zeta function

In 2001 Boudjelkha [12] also developed the following Hurwitz type formula for the alternating Hurwitz zeta function $\varsigma_a(s, t)$ which he defined for $\sigma > 0$, $0 < t \le 1$ by

(6.1) $$\varsigma_a(s,t) = \frac{1}{\Gamma(s)} \int_0^\infty \frac{x^{s-1}e^{-x(t-1)}}{e^x + 1} dx$$

where, as usual, $\sigma = \mathrm{Re}(s)$. When $t = 1$ we have for $\sigma > 0$, $\varsigma_a(s,1) = \varsigma_a(s)$ [36, p.103].

Boudjelkha's formula is

(6.2) $$\varsigma_a(s,t) = 2\Gamma(1-s)\pi^{s-1} \left[ \sin\left(\frac{\pi s}{2}\right) \sum_{n=0}^\infty \frac{\cos(2n+1)\pi t}{(2n+1)^{1-s}} + \cos\left(\frac{\pi s}{2}\right) \sum_{n=0}^\infty \frac{\sin(2n+1)\pi t}{(2n+1)^{1-s}} \right]$$

and holds under the same conditions as (1.3) above, namely:

(6.3) $$(\sigma < 0, \, 0 < t \le 1; \; 0 < \sigma, \, t < 1)$$

This may be written as

(6.4) $$\varsigma_a(s,t) = 2\Gamma(1-s)\pi^{s-1} \sum_{n=0}^\infty \frac{\cos[\pi s/2 - (2n+1)\pi t]}{(2n+1)^{1-s}}$$

The notation $\eta(s,t)$ is frequently used instead of $\varsigma_a(s,t)$.

Guillera and Sondow [24] proved that for all complex values of $s$ and complex $z$ such that $\mathrm{Re}(z) < \frac{1}{2}$

(6.5) $$(1-z)\Phi(z,s,t) = \sum_{n=0}^\infty \left(\frac{-z}{1-z}\right)^n \sum_{k=0}^n \binom{n}{k} \frac{(-1)^k}{(t+k)^s}$$

where $\Phi(z,s,t)$ the Hurwitz-Lerch zeta function $\Phi(z,s,x)$ defined by [36, p.121] as

(6.6) $$\Phi(z,s,t) = \sum_{n=0}^\infty \frac{z^n}{(n+t)^s}$$

With $z = -1$ we obtain

(6.7) $$\Phi(-1,s,t) = \sum_{n=0}^\infty \frac{1}{2^{n+1}} \sum_{k=0}^n \binom{n}{k} \frac{(-1)^k}{(t+k)^s}$$

which was determined in a different manner in equation (4.4.79) in [16].

With $t = 1$ we have

$$(6.8) \qquad \sum_{n=0}^{\infty} \frac{1}{2^{n+1}} \sum_{k=0}^{n} \binom{n}{k} \frac{(-1)^k}{(k+1)^s} = \sum_{n=0}^{\infty} \frac{(-1)^n}{(n+1)^s} = \sum_{n=1}^{\infty} \frac{(-1)^{n+1}}{n^s} = \varsigma_a(s)$$

which is the Hasse/Sondow identity (see [27] and [34]).

We note that

$$\frac{e^{-x(t-1)}}{e^x - z} = \sum_{n=0}^{\infty} z^n e^{-(t+n)x}$$

and

$$\int_0^{\infty} e^{-(t+n)y} x^{s-1} dx = \frac{1}{(t+n)^s} \int_0^{\infty} e^{-u} u^{s-1} du = \frac{\Gamma(s)}{(t+n)^s}$$

and we therefore have the integral representation [36, p.121]

$$(6.9) \qquad \Phi(z,s,t) = \frac{1}{\Gamma(s)} \int_0^{\infty} \frac{x^{s-1} e^{-x(t-1)}}{e^x - z} dx$$

Comparing (6.1), (6.5) and (6.9) we deduce that

$$(6.10) \qquad \Phi(-1,s,t) = \sum_{n=0}^{\infty} \frac{(-1)^n}{(n+t)^s} = \varsigma_a(s,t) = \sum_{n=0}^{\infty} \frac{1}{2^{n+1}} \sum_{k=0}^{n} \binom{n}{k} \frac{(-1)^k}{(t+k)^s}$$

It appears that this formula arises as the result of applying the Euler series transformation [28, p.244].

Williams and Zhang [40] also considered the alternating Hurwitz zeta function in 1993 (but their paper was not referenced by Boudjelkha [12]).Williams and Zhang defined $J(s,t)$ by

$$J(s,t) = \sum_{n=0}^{\infty} \frac{(-1)^n}{(n+t)^s}$$

and they reported that for $\sigma < 0$

$$J(s,t) = 2\Gamma(1-s)\pi^{s-1} \left[ \sin\left(\frac{\pi s}{2}\right) \sum_{n=0}^{\infty} \frac{\cos(2n+1)\pi t}{(2n+1)^{1-s}} + \cos\left(\frac{\pi s}{2}\right) \sum_{n=0}^{\infty} \frac{\sin(2n+1)\pi t}{(2n+1)^{1-s}} \right]$$

(which is here reproduced after inserting a factor of 2 which appears to be missing in equations (1.7) and (3.4) of their paper [40]). There appears to be some confusion in equation (1.7) of [40] which states that it is valid for $\sigma < 1$ whereas equation (3.3) of the same paper states that the requisite condition is $\sigma < 0$.

Upon a separation of terms according to the parity of $n$ we see that for $\sigma > 1$

$$\sum_{n=0}^{\infty} \frac{(-1)^n}{(n+t)^s} = \sum_{n=0}^{\infty} \frac{1}{(2n+t)^s} - \sum_{n=0}^{\infty} \frac{1}{(2n+1+t)^s}$$

$$= 2^{-s} \left[ \sum_{n=0}^{\infty} \frac{1}{(n+t/2)^s} - \sum_{n=0}^{\infty} \frac{1}{(n+(t+1)/2)^s} \right]$$

and we therefore see that $\varsigma_a(s,t)$ is related to the Hurwitz zeta function by the formula

(6.11)  $$\varsigma_a(s,t) = 2^{-s} \left[ \varsigma\left(s, \frac{t}{2}\right) - \varsigma\left(s, \frac{1+t}{2}\right) \right]$$

Hansen and Patrick [25] showed in 1962 that the Hurwitz zeta function could be written as

(6.12)  $$\varsigma(s,x) = 2^s \varsigma(s,2x) - \varsigma\left(s, x+\frac{1}{2}\right)$$

and, by analytic continuation, this holds for all $s$. With $x = t/2$ this becomes

(6.13)  $$\varsigma\left(s, \frac{t}{2}\right) = 2^s \varsigma(s,t) - \varsigma\left(s, \frac{1+t}{2}\right)$$

and hence we have for $\sigma > 1$

(6.14)  $$\varsigma_a(s,t) = \varsigma(s,t) - 2^{1-s} \varsigma\left(s, \frac{1+t}{2}\right)$$

and

(6.15)  $$\varsigma_a(s,t) = 2^{1-s} \varsigma\left(s, \frac{t}{2}\right) - \varsigma(s,t)$$

Since $\varsigma(s,t)$ can be continued analytically to the whole complex plane except for a simple pole at $s = 1$, $\varsigma_a(s,t)$ can be continued analytically to become an entire function and (6.11), (6.14) and (6.15) therefore hold in the whole complex plane.

We now multiply (c) by $s-1$

$$(s-1)\varsigma_a(s,t) = 2^{1-s}(s-1)\varsigma\left(s, \frac{t}{2}\right) - (s-1)\varsigma(s,t)$$

and take the limit as $s \to 1$ to obtain

(6.16) $$\lim_{s \to 1}[(s-1)\varsigma_a(s,t)] = 0$$

since $\lim_{s \to 1}[(s-1)\varsigma(s,t)] = 1$.

It should be noted that we cannot automatically substitute $s = 0$ in the formula $\varsigma_a(s) = (1-2^{1-s})\varsigma(s)$ because that equation is only valid for $\operatorname{Re}(s) > 0$ (excluding $s = 1$). Fortunately, Hardy [38, p.16] gave the following functional equation for the alternating zeta function

(6.17) $$\varsigma_a(-s) = \left(1 - \left[2^{-s} - 1\right]^{-1}\right)\pi^{-s-1}s\Gamma(s)\sin(\pi s/2)\varsigma_a(1+s)$$

$$= 2\frac{\left[2^{-s-1} - 1\right]}{\left[2^{-s} - 1\right]}\pi^{-s-1}s\Gamma(s)\sin(\pi s/2)\varsigma_a(1+s)$$

and it is this equation that enables us to equate $\varsigma_a(0) = -\varsigma(0)$. As can be seen from Ayoub's paper [6], this is precisely the functional equation for the zeta function which was first postulated by Euler many years before Riemann.

We have

$$\varsigma_a(0) = -\pi^{-1}\varsigma_a(1)\lim_{s \to 0}\frac{\sin(\pi s/2)}{\left[2^{-s} - 1\right]}$$

Using L'Hôpital's rule results in

$$= \pi^{-1}\varsigma_a(1)\frac{\pi}{2}\lim_{s \to 0}\frac{\cos(\pi s/2)}{2^{-s}\log 2} = \frac{\varsigma_a(1)}{2\log 2}$$

We note that

$$\varsigma_a(1) = \lim_{s \to 1}[(1-2^{1-s})\varsigma(s)]$$

$$= \lim_{s \to 1}\left[\frac{1-2^{1-s}}{s-1}(s-1)\varsigma(s)\right]$$

$$= \lim_{s \to 1}\left[\frac{1-2^{1-s}}{s-1}\right]\lim_{s \to 1}[(s-1)\varsigma(s)]$$

and using L'Hôpital's rule again gives us

$$= \log 2$$

We then have the well-known result

$$\varsigma_a(0) = \frac{1}{2}$$

It is interesting to note that substituting (1.3) in (6.11) gives us

$$\varsigma_a(s,t) = 2^{1-s}\Gamma(1-s)\left[\sin\left(\frac{\pi s}{2}\right)\sum_{n=1}^{\infty}\frac{\cos n\pi t - \cos n\pi(1+t)}{(2\pi n)^{1-s}} + \cos\left(\frac{\pi s}{2}\right)\sum_{n=1}^{\infty}\frac{\sin n\pi t - \sin n\pi(1+t)}{(2\pi n)^{1-s}}\right]$$

$$= 2^{1-s}\Gamma(1-s)\left[\sin\left(\frac{\pi s}{2}\right)\sum_{n=1}^{\infty}\frac{[1-(-1)^n]\cos n\pi t}{(2\pi n)^{1-s}} + \cos\left(\frac{\pi s}{2}\right)\sum_{n=1}^{\infty}\frac{[1-(-1)^n]\sin n\pi t}{(2\pi n)^{1-s}}\right]$$

$$= 2\Gamma(1-s)\pi^{1-s}\left[\sin\left(\frac{\pi s}{2}\right)\sum_{n=1}^{\infty}\frac{\cos(2n+1)\pi t}{(2n+1)^{1-s}} + \cos\left(\frac{\pi s}{2}\right)\sum_{n=1}^{\infty}\frac{\sin(2n+1)\pi t}{(2n+1)^{1-s}}\right]$$

and we have therefore recovered (6.2) in a rather straightforward manner.

Letting $s \to 1-s$ in (6.2) gives us

$$\varsigma_a(1-s,t) = 2\Gamma(s)\pi^s\left[\cos\left(\frac{\pi s}{2}\right)\sum_{n=1}^{\infty}\frac{\cos(2n+1)\pi t}{(2n+1)^s} + \sin\left(\frac{\pi s}{2}\right)\sum_{n=1}^{\infty}\frac{\sin(2n+1)\pi t}{(2n+1)^s}\right]$$

and using (6.10) this is equal to

$$= \sum_{n=0}^{\infty}\frac{1}{2^{n+1}}\sum_{k=0}^{n}\binom{n}{k}(-1)^k(t+k)^{s-1}$$

The Euler polynomials $E_m(t)$ may be expressed by

$$(6.18) \qquad E_m(t) = \sum_{n=0}^{\infty}\frac{1}{2^n}\sum_{k=0}^{n}\binom{n}{k}(-1)^k(t+k)^m$$

(which may be contrasted with equation (2.14)) and hence we obtain the well-known Fourier series [40]

(6.19a)
$$\sum_{n=1}^{\infty} \frac{\sin(2n+1)\pi t}{(2n+1)^{2m+1}} = \frac{(-1)^m \pi^{2m+1} E_{2m}(t)}{4(2m)!}$$

(6.19b)
$$\sum_{n=1}^{\infty} \frac{\cos(2n+1)\pi t}{(2n+1)^{2m}} = \frac{(-1)^m \pi^{2m} E_{2m-1}(t)}{4(2m-1)!}$$

## 7. Some trigonometric series

The following identities are recorded by Hansen [26, pp. 223 & 244] for $\mathrm{Re}(s) > 1$ and $0 < x < 2\pi$

(7.1a)

$$\sum_{n=1}^{\infty} \frac{\sin(nx+y)}{n^s} = \frac{(2\pi)^s}{2\Gamma(s)} \mathrm{cosec}(\pi s) \left[ \cos\left(y - \frac{\pi s}{2}\right) \varsigma\left(1-s, \frac{x}{2\pi}\right) - \cos\left(y + \frac{\pi s}{2}\right) \varsigma\left(1-s, 1 - \frac{x}{2\pi}\right) \right]$$

(7.1b)

$$\sum_{n=1}^{\infty} \frac{\cos(nx+y)}{n^s} = \frac{(2\pi)^s}{2\Gamma(s)} \mathrm{cosec}(\pi s) \left[ \sin\left(y + \frac{\pi s}{2}\right) \varsigma\left(1-s, 1 - \frac{x}{2\pi}\right) - \sin\left(y - \frac{\pi s}{2}\right) \varsigma\left(1-s, \frac{x}{2\pi}\right) \right]$$

and we may note that the second identity may be obtained by differentiating the first one with respect to $y$.

Note that as $s \to 2N$ in the first identity, $\mathrm{cosec}(\pi s) \to \infty$; this point is considered in more detail later.

Letting $y = 0$ and $x \to \pi x$ we obtain

(7.2a) $\quad \sum_{n=1}^{\infty} \frac{\sin(n\pi x)}{n^s} = \frac{(2\pi)^s}{2\Gamma(s)} \mathrm{cosec}(\pi s) \cos\left(\frac{\pi s}{2}\right) \left[ \varsigma\left(1-s, \frac{x}{2}\right) - \varsigma\left(1-s, 1-\frac{x}{2}\right) \right]$

(7.2b) $\quad \sum_{n=1}^{\infty} \frac{\cos(n\pi x)}{n^s} = \frac{(2\pi)^s}{2\Gamma(s)} \mathrm{cosec}(\pi s) \sin\left(\frac{\pi s}{2}\right) \left[ \varsigma\left(1-s, 1-\frac{x}{2}\right) + \varsigma\left(1-s, \frac{x}{2}\right) \right]$

or equivalently as

(7.3a) $\quad \sum_{n=1}^{\infty} \frac{\sin(n\pi x)}{n^s} = \frac{(2\pi)^s}{4\Gamma(s)} \mathrm{cosec}\left(\frac{\pi s}{2}\right) \left[ \varsigma\left(1-s, \frac{x}{2}\right) - \varsigma\left(1-s, 1-\frac{x}{2}\right) \right]$

(7.3b) $\qquad \sum_{n=1}^{\infty} \frac{\cos(n\pi x)}{n^s} = \frac{(2\pi)^s}{4\Gamma(s)} \sec\left(\frac{\pi s}{2}\right)\left[\varsigma\left(1-s, 1-\frac{x}{2}\right) + \varsigma\left(1-s, \frac{x}{2}\right)\right]$

Letting $x \to 2t$ and $p = 1 - s$ and then adding the above two equations immediately results in the well-known Hurwitz's formula for the Fourier expansion of the Riemann zeta function $\varsigma(p, t)$

(7.4) $\qquad \varsigma(p,t) = 2\Gamma(1-p)\left[\sin\left(\frac{\pi p}{2}\right)\sum_{n=1}^{\infty}\frac{\cos 2n\pi t}{(2\pi n)^{1-p}} + \cos\left(\frac{\pi p}{2}\right)\sum_{n=1}^{\infty}\frac{\sin 2n\pi t}{(2\pi n)^{1-p}}\right]$

where $\text{Re}(p) < 0$ and $0 < t \leq 1$. Boudjelkha [12] showed that this formula also applies in the region $\text{Re}(p) < 1$. It may be noted that when $t = 1$ this reduces to Riemann's functional equation for $\varsigma(p)$. Letting $s = 1 - p$ we may write this as

(7.5) $\qquad \varsigma(1-s,t) = 2\Gamma(s)\left[\cos\left(\frac{\pi s}{2}\right)\sum_{n=1}^{\infty}\frac{\cos 2n\pi t}{(2\pi n)^s} + \sin\left(\frac{\pi s}{2}\right)\sum_{n=1}^{\infty}\frac{\sin 2n\pi t}{(2\pi n)^s}\right]$

This may also be written as [36, p.89]

(7.6) $\qquad \varsigma(1-s,t) = \frac{\Gamma(s)}{(2\pi)^s}\left[e^{-\frac{1}{2}\pi i t}L(t,s) + e^{\frac{1}{2}\pi i t}L(-t,s)\right]$

where the periodic (or Lerch) zeta function $L(t,s)$ is defined by

$$L(t,s) = \sum_{n=1}^{\infty}\frac{e^{2\pi i t}}{n^s}$$

We now consider the equations (7.3) in the case where $s$ is a positive integer $N > 1$: we see that

(7.7a) $\qquad \sum_{n=1}^{\infty} \frac{\sin(n\pi x)}{n^N} = \frac{(2\pi)^N}{4(N-1)!}\operatorname{cosec}\left(\frac{\pi N}{2}\right)\left[\varsigma\left(1-N, \frac{x}{2}\right) - \varsigma\left(1-N, 1-\frac{x}{2}\right)\right]$

(7.7b) $\qquad \sum_{n=1}^{\infty} \frac{\cos(n\pi x)}{n^N} = \frac{(2\pi)^N}{4(N-1)!}\sec\left(\frac{\pi N}{2}\right)\left[\varsigma\left(1-N, 1-\frac{x}{2}\right) + \varsigma\left(1-N, \frac{x}{2}\right)\right]$

Using the familiar identity [5, p.264] for $N \geq 1$ (which is also derived in (7.12) below)

$$\varsigma(1-N,t) = -\frac{B_N(t)}{N}$$

we therefore have

$$\varsigma\left(1-N,\frac{x}{2}\right)-\varsigma\left(1-N,1-\frac{x}{2}\right)=-\frac{1}{N}\left[B_N\left(\frac{x}{2}\right)-B_N\left(1-\frac{x}{2}\right)\right]$$

In passing we note that

$$\lim_{N\to 0}N\varsigma(1-N,t)=-1=-\lim_{N\to 0}B_N(t)$$

Using the well-known formula for the Bernoulli polynomials [36, p.60]

$$B_n(1-t)=(-1)^n B_n(t)$$

we get for $N\geq 1$

$$\varsigma\left(1-N,\frac{x}{2}\right)-\varsigma\left(1-N,1-\frac{x}{2}\right)=\frac{1}{N}\left[(-1)^N-1\right]B_N\left(\frac{x}{2}\right)$$

In a similar manner we see that for $N\geq 1$

$$\varsigma\left(1-N,1-\frac{x}{2}\right)+\varsigma\left(1-N,\frac{x}{2}\right)=\frac{1}{N}\left[(-1)^{N+1}-1\right]B_N\left(\frac{x}{2}\right)$$

Accordingly the above equations (7.7) may be written as

(7.8a) $$\sum_{n=1}^{\infty}\frac{\sin(n\pi x)}{n^N}=\frac{(2\pi)^N}{4N!}\operatorname{cosec}\left(\frac{\pi N}{2}\right)\left[(-1)^N-1\right]B_N\left(\frac{x}{2}\right)$$

(7.8b) $$\sum_{n=1}^{\infty}\frac{\cos(n\pi x)}{n^N}=\frac{(2\pi)^N}{4N!}\sec\left(\frac{\pi N}{2}\right)\left[(-1)^{N+1}-1\right]B_N\left(\frac{x}{2}\right)$$

Letting $N\to 2N+1$ and $x=2t$ in (7.8a) then gives us the well-known Fourier series for the odd Bernoulli polynomials [4, p.338]

(7.9a) $$B_{2N+1}(t)=(-1)^{N+1}2(2N+1)!\sum_{n=1}^{\infty}\frac{\sin 2n\pi t}{(2\pi n)^{2N+1}}\qquad,N=0,1,2,...$$

and letting $N\to 2N$ and $x=2t$ in (7.8b) gives us the corresponding series for the even Bernoulli polynomials

(7.9b) $$B_{2N}(t)=(-1)^{N+1}2(2N)!\sum_{n=1}^{\infty}\frac{\cos 2n\pi t}{(2\pi n)^{2N}}\qquad,N=1,2,...$$

Referring to (7.7a) we see that we have an indeterminate form when $N \rightarrow 2N$ and hence we step back to (7.3a)

$$\sum_{n=1}^{\infty} \frac{\sin(n\pi x)}{n^s} = \frac{(2\pi)^s}{4\Gamma(s)} \operatorname{cosec}\left(\frac{\pi s}{2}\right)\left[\varsigma\left(1-s,\frac{x}{2}\right) - \varsigma\left(1-s,1-\frac{x}{2}\right)\right]$$

and consider the limit as $s \rightarrow 2N$. Using L'Hôpital's rule we see that

$$\lim_{s \rightarrow 2N} \operatorname{cosec}\left(\frac{\pi s}{2}\right)\left[\varsigma\left(1-s,\frac{x}{2}\right) - \varsigma\left(1-s,1-\frac{x}{2}\right)\right] = \lim_{s \rightarrow 2N} \frac{\left[-\varsigma'\left(1-s,\frac{x}{2}\right) + \varsigma'\left(1-s,1-\frac{x}{2}\right)\right]}{\frac{\pi}{2}\cos\left(\frac{\pi s}{2}\right)}$$

$$= \frac{2}{\pi}\left[-\varsigma'\left(1-2N,\frac{x}{2}\right) + \varsigma'\left(1-2N,1-\frac{x}{2}\right)\right]$$

and we therefore obtain

(7.10a) $$\sum_{n=1}^{\infty} \frac{\sin(n\pi x)}{n^{2N}} = \frac{(2\pi)^{2N-1}}{(2N-1)!}\left[\varsigma'\left(1-2N,1-\frac{x}{2}\right) - \varsigma'\left(1-2N,\frac{x}{2}\right)\right]$$

Applying the same limiting procedure to (7.3b) gives us

(7.10b) $$\sum_{n=1}^{\infty} \frac{\cos(n\pi x)}{n^{2N+1}} = (-1)^N \frac{(2\pi)^{2N}}{(2N)!}\left[\varsigma'\left(-2N,1-\frac{x}{2}\right) + \varsigma'\left(-2N,\frac{x}{2}\right)\right]$$

With $N = 0$ in (7.10b) we obtain

$$\sum_{n=1}^{\infty} \frac{\cos(n\pi x)}{n} = \varsigma'\left(0,1-\frac{x}{2}\right) + \varsigma'\left(0,\frac{x}{2}\right)$$

and using Lerch's identity (2.16) this becomes

$$= \log\Gamma\left(1-\frac{x}{2}\right) + \log\Gamma\left(\frac{x}{2}\right) - \log(2\pi)$$

Employing Euler's reflection formula we obtain

$$= -\log\sin\left(\frac{\pi x}{2}\right) - \log 2$$

and we thereby obtain the familiar Fourier series [13, p.241]

$$\log\left[2\sin\left(\frac{\pi x}{2}\right)\right] = -\sum_{n=1}^{\infty}\frac{\cos(n\pi x)}{n}$$

With $N = 1$ in (7.8a) we see that

$$\frac{\pi}{2}(1-x) = \sum_{n=1}^{\infty}\frac{\sin(n\pi x)}{n}$$

The Clausen functions $\text{Cl}_N(x)$ are defined by [36, p.115]

$$\text{Cl}_{2N}(x) = \sum_{n=1}^{\infty}\frac{\sin nx}{n^{2N}}$$

$$\text{Cl}_{2N+1}(x) = \sum_{n=1}^{\infty}\frac{\cos nx}{n^{2N+1}}$$

and using (7.10) we have thereby derived Adamchik's results [2]

(7.11a) $\qquad \dfrac{(2N-1)!}{(2\pi)^{2N-1}}\text{Cl}_{2N}(2\pi x) = \varsigma'\left(1-2N,x\right) - \varsigma'\left(1-2N,1-x\right)$

(7.11b) $\qquad (-1)^N\dfrac{(2N)!}{(2\pi)^{2N}}\text{Cl}_{2N+1}(2\pi x) = \varsigma'\left(-2N,x\right) + \varsigma'\left(-2N,1-x\right)$

which were also obtained in equation (4.3.167) in [14] by an entirely different method.

We recall Hasse's formula [27] for the Hurwitz zeta function which is valid for all $s \in \mathbf{C}$ except for $s = 1$

$$(s-1)\varsigma(s,u) = \sum_{n=0}^{\infty}\frac{1}{n+1}\sum_{k=0}^{n}\binom{n}{k}\frac{(-1)^k}{(u+k)^{s-1}}$$

and with $s \to 1-s$ this may be written as

$$-s\varsigma(1-s,u) = \sum_{n=0}^{\infty}\frac{1}{n+1}\sum_{k=0}^{n}\binom{n}{k}(-1)^k(u+k)^s$$

It is shown in equation (A.23) of [19] that

$$B_m(x) = \sum_{n=0}^{m}\frac{1}{n+1}\sum_{k=0}^{n}(-1)^k\binom{n}{k}(x+k)^m$$

A different proof, using the Hurwitz-Lerch zeta function, was recently given by Guillera and Sondow [24]. They also noted that

$$\sum_{k=0}^{n}(-1)^k\binom{n}{k}(x+k)^m=0 \quad \text{for } n> m=0,1,2,...$$

and we therefore have

$$B_m(x)=\sum_{n=0}^{\infty}\frac{1}{n+1}\sum_{k=0}^{n}(-1)^k\binom{n}{k}(x+k)^m$$

Using the Hasse identity this immediately gives us the well-known result

(7.12) $$\varsigma(1-m,x)=-\frac{B_m(x)}{m}$$

which we used above.

With $x=\dfrac{2p}{q}$ in (7.3a) we get

$$\sum_{n=1}^{\infty}\frac{\sin(2n\pi p/q)}{n^s}=\frac{(2\pi)^s}{4\Gamma(s)}\operatorname{cosec}\left(\frac{\pi s}{2}\right)\left[\varsigma\left(1-s,\frac{p}{q}\right)-\varsigma\left(1-s,1-\frac{p}{q}\right)\right]$$

and from (2.1.1) we have

$$\varsigma\left(s,\frac{p}{q}\right)-\varsigma\left(s,1-\frac{p}{q}\right)=4\Gamma(1-s)(2\pi q)^{s-1}\cos\left(\frac{\pi s}{2}\right)\sum_{j=1}^{q}\sin\left(\frac{2\pi jp}{q}\right)\varsigma\left(1-s,\frac{j}{q}\right)$$

$$\varsigma\left(1-s,\frac{p}{q}\right)-\varsigma\left(1-s,1-\frac{p}{q}\right)=4\Gamma(s)(2\pi q)^{-s}\sin\left(\frac{\pi s}{2}\right)\sum_{j=1}^{q}\sin\left(\frac{2\pi jp}{q}\right)\varsigma\left(s,\frac{j}{q}\right)$$

We then see that

$$\sum_{n=1}^{\infty}\frac{\sin(2n\pi p/q)}{n^s}=q^{-s}\sum_{j=1}^{q}\sin\left(\frac{2\pi jp}{q}\right)\varsigma\left(s,\frac{j}{q}\right)$$

and this generalises the formula previously given by  Srivastava and Tsumura [37].

Letting $x=2\pi t$ and $y=\pi t$ in (7.1a) gives us

$$\sum_{n=1}^{\infty} \frac{\sin(2n+1)\pi t}{n^s} = \frac{(2\pi)^s}{2\Gamma(s)} \operatorname{cosec}(\pi s)\left[\cos\left(\pi t - \frac{\pi s}{2}\right)\varsigma(1-s,t) - \cos\left(\pi t + \frac{\pi s}{2}\right)\varsigma(1-s,1-t)\right]$$

## 9. Some connections with the sine and cosine integrals

In passing, we mention that there exist other trigonometric expansions for $\psi(x), \log\Gamma(x)$, $\log G(1+x)$, $\varsigma'(-1,x)$ etc. These are set out below (further details are contained in [18]).

(8.1) $$\psi(x) = \log x - \frac{1}{2x} + 2\sum_{n=1}^{\infty}[\cos(2n\pi x)Ci(2n\pi x) + \sin(2n\pi x)si(2n\pi x)]$$

which appears in Nörlund's book [33, p.108].

(8.2) $$\log\Gamma(x) =$$

$$\frac{1}{2}\log(2\pi) + \left(x - \frac{1}{2}\right)\log x - x + \frac{1}{\pi}\sum_{n=1}^{\infty}\frac{1}{n}[\sin(2n\pi x)Ci(2n\pi x) - \cos(2n\pi x)si(2n\pi x)]$$

(8.3) $$x\log\Gamma(x) - \log G(1+x) = \frac{1}{4}x\left[-x + 2(x-1)\log x\right] + \frac{\text{Cl}_2(2\pi x)}{4\pi}$$

$$-\frac{1}{2\pi^2}\sum_{n=1}^{\infty}\frac{1}{n^2}\left[\cos(2n\pi x)Ci(2n\pi x) + \sin(2n\pi x)Si(2n\pi x)\right] + \frac{1}{12} - \varsigma'(-1)$$

Elizalde [20] reported in 1985 that for $x > 0$

(8.4) $$\varsigma'(-1,x) =$$

$$-\varsigma(-1,x)\log x - \frac{1}{4}x^2 + \frac{1}{12} - \frac{1}{2\pi^2}\sum_{n=1}^{\infty}\frac{1}{n^2}[\cos(2n\pi x)Ci(2n\pi x) + \sin(2n\pi x)si(2n\pi x)]$$

where $Si(x)$ is the sine integral function defined by [23, p.878] and [1, p.231] as

$$Si(x) = \int_0^x \frac{\sin t}{t}\,dt \qquad , \quad Si(0) = 0$$

We have the well-known integral from Fourier series analysis

$$\frac{\pi}{2} = \int_0^\infty \frac{\sin t}{t}\,dt$$

and therefore defining

$$si(x) = Si(x) - \frac{\pi}{2}$$

we have

$$si(x) = \int_0^x \frac{\sin t}{t}\,dt - \int_0^\infty \frac{\sin t}{t}\,dt = -\int_x^\infty \frac{\sin t}{t}\,dt$$

The cosine integral $Ci(x)$ is defined in [23, p.878] and [1, p.231] as

$$Ci(x) = \gamma + \log x + \int_0^x \frac{\cos t - 1}{t}\,dt = \gamma + \log x + \sum_{n=1}^\infty \frac{(-1)^n x^{2n}}{2n(2n)!}$$

(where, in the final part, we have simply substituted the Maclaurin series for the integrand). We also have

$$Ci(x) = -\int_x^\infty \frac{\cos t}{t}\,dt$$

and this more clearly shows the connection with $si(x)$. $Ci(x)$ is frequently designated as $ci(x)$ in other works such as [23].

$\square$

We recall the well-known identity

(8.5)
$$\gamma = \lim_{n\to\infty} \sum_{k=1}^n \left[ \frac{1}{k} - \log\left(1 + \frac{1}{k}\right) \right]$$

It is interesting to note that Sondow [35] has discovered a similar alternating series

(8.6)
$$\log\frac{4}{\pi} = \lim_{n\to\infty} \sum_{k=1}^n (-1)^{k+1} \left[ \frac{1}{k} - \log\left(1 + \frac{1}{k}\right) \right]$$

$$= \sum_{k=1}^\infty \frac{(-1)^{k+1}}{k} - \sum_{k=1}^\infty (-1)^{k+1} \log\left(1 + \frac{1}{k}\right)$$

$$= \log 2 - \log \frac{\pi}{2} = \log \frac{4}{\pi}$$

where, in the last line, we have employed (4.4.102) from [16].

Sondow's formula may also be obtained from Lerch's trigonometric series expansion for the digamma function for $0 < x < 1$ (see for example Gronwall's paper [22, p.105] and Nielsen's book [32, p.204])

(8.7) $\qquad \psi(x)\sin \pi x + \frac{\pi}{2}\cos \pi x + (\gamma + \log 2\pi)\sin \pi x = \sum_{n=1}^{\infty} \sin(2n+1)\pi x . \log \frac{n+1}{n}$

Letting $x = 1/2$ we obtain

$$\psi(1/2) + \gamma + \log 2\pi = \sum_{n=1}^{\infty} (-1)^n \log \frac{n+1}{n}$$

and, since [36, p.20] $\psi(1/2) = -\gamma - 2\log 2$, equation (8.6) follows automatically.

We may also write (8.7) as

(8.8) $\qquad \psi(x) + \frac{\pi}{2}\cot \pi x + (\gamma + \log 2\pi) = \sum_{n=1}^{\infty} \frac{\sin(2n+1)\pi x}{\sin \pi x} . \log \frac{n+1}{n}$

The above formula suggests that integration may be fruitful employing the integral in the tables [23, p.163]

$$\int \frac{\sin(2n+1)x}{\sin x} dx = x + \sum_{k=1}^{n} \frac{\sin 2kx}{k}$$

The integral $\int x^p \psi(x)$ may also produce interesting results.

One could also employ the substitution in (8.7)

$$\log \frac{n+1}{n} = \int_0^1 \frac{y^n(y-1)}{y \log y} dy$$

but this has not been explored in any depth.

Introduction to the General Theory of Infinite Processes and of Analytic Functions; With an Account of the Principal Transcendental Functions. Fourth Ed., Cambridge University Press, Cambridge, London and New York, 1963.

Donal F. Connon
Elmhurst
Dundle Road
Matfield
Kent TN12 7HD
dconnon@btopenworld.com